\documentstyle[amssymb,amsfonts]{amsart}

\def\R{{\hbox{\bf R}}}

 at 10 true pt

\def\allt#1{%
\smash{
 \vtop{%
     \ialign{%
        ##\crcr
        $\hfil\displaystyle{\tilde \forall}\hfil$\crcr%
        \noalign{\kern1.5pt\nointerlineskip}
        $\hfil\!\!#1\hfil$\crcr\noalign{\kern1.5pt}
        }
       }
      } \hbox{$\vphantom{#1}$}
     }

\def\be#1{\begin{equation}\label{#1}}
\def\bas{\begin{align*}}
\def\eas{\end{align*}}
\def\bi{\begin{itemize}}
\def\ei{\end{itemize}}

\def\Z{{\hbox{\bf Z}}}

\def\eps{\varepsilon}
\newenvironment{proof}{\noindent {\bf Proof} }{\endprf\par}
\def \endprf{\hfill  {\vrule height6pt width6pt depth0pt}\medskip}
\def\emph#1{{\it #1}}
\def\textbf#1{{\bf #1}}

\theoremstyle{plain}
  \newtheorem{theorem}[subsection]{Theorem}

  \newtheorem{proposition}[subsection]{Proposition}
  
  \newtheorem{lemma}[subsection]{Lemma}
  \newtheorem{corollary}[subsection]{Corollary}

\theoremstyle{remark}

\theoremstyle{definition}

\include{psfig}

\begin{document}

\title[Sum-product estimate]{A sum-product estimate in finite fields, and applications}

\author{Jean Bourgain}
\address{School of Mathematics, Institute of Advanced Study, Princeton NJ 08540}
\email{bourgain@@math.ias.edu}

\author{Nets Katz}
\address{Department of Mathematics, Washington University in St. Louis, St. Louis MO 63130}
\email{nets@@math.wustl.edu}

\author{Terence Tao}
\address{Department of Mathematics, UCLA, Los Angeles CA 90095-1555}
\email{tao@@math.ucla.edu}

\begin{abstract}
Let $A$ be a subset of a finite field $F := \Z/q\Z$ for some prime $q$.  If $|F|^\delta < |A| < |F|^{1-\delta}$ for some $\delta > 0$, then we prove the estimate $|A+A| + |A \cdot A| \geq c(\delta) |A|^{1+\eps}$ for some $\eps = \eps(\delta) > 0$.  This is a finite field analogue of a result of \cite{esz}.  We then use this estimate to prove a Szemer\'edi-Trotter type theorem in finite fields, and obtain a new estimate for the Erd\"os distance problem in finite fields, as well as the three-dimensional Kakeya problem in finite fields.
\end{abstract}
\maketitle

\section{Introduction}

Let $A$ be a non-empty subset of a finite field $F$.  We consider the sum set
$$ A+A := \{ a+b: a, b \in A \}$$
and the product set
$$ A \cdot A := \{ a \cdot b: a, b \in A \}.$$
Let $|A|$ denote the cardinality of $A$.  Clearly we have the bounds
$$ |A+A|, |A \cdot A| \geq |A|.$$
These bounds are clearly sharp when $A$ is a subfield of $F$; however when $A$ is not a subfield (or an affine transformation of a subfield) then we expect some improvement.  In particular, when $F$ is the cyclic field $F := \Z / q\Z$ for some prime $q$, then $F$ has no proper subfields, and one expects some gain when $1 \ll |A| \ll |F|$.  The first main result of this paper is to show that this is indeed the case:

\begin{theorem}[Sum-product estimate]\label{sum-product}  Let $F := \Z / q\Z$ for some prime $q$, and let $A$ be a subset of $F$ such that
$$ |F|^\delta < |A| < |F|^{1-\delta}$$
for some $\delta > 0$.  Then one has a bound of the form
\begin{equation}\label{sp}
\max(|A+A|, |A \cdot A|) \geq c(\delta) |A|^{1+\eps}
\end{equation}
for some $\eps = \eps(\delta) > 0$.
\end{theorem}

We note that one needs both $|A+A|$ and $|A \cdot A|$ on the left-hand side to obtain an estimate of this type; for the additive term $|A+A|$ this can be seen by considering an arithmetic progression such as $A := \{1, \ldots, N \}$, and for the multiplicative term $|A \cdot A|$ this can be seen by considering a geometric progression.  Thus the above estimate can be viewed as a statement that a set cannot behave like an arithmetic progression and a geometric progression simultaneously.  This suggests using Freiman's theorem \cite{freiman} to obtain the estimate \eqref{sp}, but the best known quantitative bounds for Freiman's theorem \cite{chang} are only able to gain a logarithmic factor in $|A|$ over the trivial bound, as opposed to the polynomial gain of $|A|^\eps$ in our result.

We do not know what the optimal value of $\eps$ should be.  If the finite field $F$ were replaced with the integers $\Z$, then it is an old conjecture of Erd\"os that one indeed has $\max(|A+A|, |A \cdot A|) \geq c(\eps) |A|^{2-\eps}$ for any $\eps > 0$, and in analogy with this integer problem one might conjecture that $\max(|A+A|, |A \cdot A|) \geq c(\eps) \min(|A|^{2-\eps}, |F|^{1-\eps})$ for all subsets $A$ of $F$.   However such an estimate, if true, is likely to be extremely difficult.  In the integer problem, the analogue of \eqref{sp} was obtained by Erd\"os and Sz\'emeredi \cite{esz}, with improvements in the value of $\eps$ by various authors; at present the best known result is $\eps = 1/4$, obtained by Elekes \cite{elekes}.  Also, a continuous version of \eqref{sp}, for fractal subsets of the real line, was recently obtained by the first author \cite{bourgain:ring}.

The proof of Theorem \ref{sum-product} is based on a recent argument of Edgar and Miller \cite{edgar}, who solved the \emph{Erd\"os ring problem} \cite{Erdos:ring}.  (An alternate solution to this problem has also appeared in \cite{bourgain:ring}).  Specifically, these authors showed that there was no Borel subring $A$ of the reals which has Hausdorff dimension strictly between 0 and 1.  Since such subrings clearly obey the identities $A+A = A \cdot A = A$, one can see that this problem has some similarities with Theorem \ref{sum-product}.  We prove this theorem in Section \ref{sp-sec}.

It has been known for some time that sum-product estimates have application to certain geometric combinatorics problems, such as the incidence problem for lines and the Erd\"os distance problem.  (See e.g. \cite{elekes}, \cite{falconer}, \cite{cst}, \cite{katz:falc}, \cite{katz:falc-expos}).  Using these ideas (and particularly those from \cite{katz:falc}, \cite{katz:falc-expos}), we can prove a theorem of Szemer\'edi-Trotter type in two-dimensional finite field geometries.  The precise statement is in Theorem \ref{sz}; roughly speaking, this theorem asserts that if we are in the finite plane $(\Z/q\Z)^2$ and one has $N$ lines and $N$ points in that plane for some $1 \ll N \ll q^2$, then there are at most $O(N^{3/2-\eps})$ incidences; this improves upon the standard bound of $O(N^{3/2})$ obtained from extremal graph theory.  We state and prove this theorem in Section \ref{sz-sec}.  Roughly speaking, the idea is to assume for contradiction that one can attain close to $N^{3/2}$ incidences, and then show that this forces most of the $N$ points to lie in a (projective transformation of a) $N^{1/2} \times N^{1/2}$ grid.  One then applies Theorem \ref{sum-product} to obtain a contradiction.  Our arguments closely follow those in \cite{katz:falc}, \cite{katz:falc-expos}.

Using this Szemer\'edi-Trotter type theorem we can also obtain a non-trivial result on the finite field Erd\"os distance problem in the case when $-1$ is not a square: specifically, we show that $N$ points in the finite plane $(\Z/q\Z)^2$ determine at least $N^{1/2+\eps}$ distinct distances if $1 \ll N \ll q^2$; this improves upon the bound of $N^{1/2}$ obtainable by extremal graph theory.  This result uses the Szem\'eredi-Trotter type theorem and the standard observation that the set of points equidistant from two fixed points $p, p'$ all lie on a line (the perpendicular bisector of $p$ and $p'$).  As such the argument is similar to those in \cite{cst}, \cite{katz:falc}, \cite{katz:falc-expos}, and is in fact rather short.  We state and prove this theorem in Section \ref{distance-sec}.

As our final application, we give a new bound on Besicovitch sets in the three-dimensional finite geometry $(\Z/q\Z)^3$.  A Besicovitch set is a set which contains a line in every direction.  The \emph{Kakeya conjecture for finite fields} asserts that such sets have cardinality at least $c(\eps) q^{3-\eps}$ for each $\eps > 0$.  The previous best lower bound known is $c q^{5/2}$, and is due to Wolff \cite{wolff:survey} (see also \cite{wolff:kakeya}, \cite{gerd:finite}).  We improve this to $c q^{5/2+\eps}$ for some absolute constant $\eps > 0$.  We prove this in Section \ref{kakeya-sec}, using some geometric ideas of the second author to transform the problem into a two-dimensional one, to which the Szemer\'edi-Trotter theorem can then be applied.
An analogous result in the continuous geometry $\R^3$ will appear by the second author elsewhere.  (An earlier result of Katz, {\L}aba, and Tao \cite{katzlabatao} also gives a similar result in the continuous case, but this result relies crucially on the fact that $\R$ has multiple scales, and so does not apply to the finite field problem).

The third author is a Clay Prize Fellow and is supported by a grant from the Packard Foundation.

\section{Some results from arithmetic combinatorics}

In this section we recall some known facts about $A+A$, $A \cdot A$, etc.  Here $F = \Z/q\Z$ is a 
finite field of prime order.

We first recall the Cauchy-Davenport inequality
\begin{equation}\label{cauchy-davenport}
|A+B|\geq \min (|A|+|B|-1, |F|)
\end{equation}
for any non-empty subsets $A,B$ of $F$.  If we are allowed to arbitrarily dilate one of the sets $A$, $B$ then we can improve
subtantially on this inequality:

\begin{lemma}\label{mixed}  Let $A$, $B$ be finite non-empty subsets of a finite field $F$, and lte $F^* := F - \{0\}$
denote the invertible elements of $F$.
Then there exists $\xi \in F^*$ such that
\begin{equation}\label{boost}
|A + B \xi| \geq \min( \frac{1}{2} |A| |B|, \frac{1}{10} |F| ).
\end{equation}
\end{lemma}

\begin{proof}
We may assume without loss of generality that $|A| |B| \leq \frac{1}{2} |F|$, since if $|A| |B| > \frac{1}{2} |F|$
we may remove some elements from $A$ and $B$ without affecting the right-hand side of
\eqref{boost}.  Let $\xi$ be an element of $F^*$.  We use the inclusion-exclusion principle\footnote{To verify our use of the principle,
suppose an element $x$ lies in $N$ of the sets $a+B\xi$ for some $N \geq 1$.  Then the sum $\sum_{a \in A} |a+B\xi|$ counts $x$ $N$ times,
while the sum $\sum_{a,a' \in A: a \neq a'} |(a+B\xi) \cap (a'+B\xi)|$ counts $x$ $N(N-1)$ times.  Since $N - \frac{N(N-1)}{2}$ is
always less than or equal to 1, the claim follows.  An alternate way to obtain this lemma (which gives slightly worse
bounds when $|A| |B| \ll |F|$, but somewhat better bounds when $|A| |B| \gg |F|$)
is by using the Cauchy-Schwarz inequality
$\| \chi_A * \chi_{B\xi}\|_{l^1}^2 \leq \| \chi_A * \chi_{B\xi}\|_{l^2} |A+B\xi|$ and again randomizing over
$\xi$.}
 and the invertibility of
$\xi$ to compute
\begin{align*}
|A + B\xi| &= |\bigcup_{a \in A} a + B\xi| \\
&\geq \sum_{a \in A} |a + B\xi| - \frac{1}{2} \sum_{a,a' \in A: a \neq a'} |(a+B\xi) \cap (a'+B\xi)| \\
&\geq \sum_{a \in A} |B| - \frac{1}{2} \sum_{a, a' \in A: a \neq a'} \sum_{b,b' \in B} \delta_{a+b\xi, a'+b'\xi}\\
&= |A| |B| - \frac{1}{2} \sum_{a, a' \in A: a \neq a'} \sum_{b,b' \in B: b \neq b'} \delta_{\xi, (a-a')/(b-b')},
\end{align*}
where $\delta_{i,j}$ is the Kronecker delta function. If we average this over all $\xi \in F_*$ we obtain
\begin{align*}
\frac{1}{|F_*|} \sum_{\xi \in F_*} |A+B\xi| &\geq |A| |B| - \frac{1}{2} \sum_{a,a' \in A: a \neq a'} \sum_{b,b' \in B: b \neq b'} \frac{1}{|F|-1}\\
&\geq |A| |B| - \frac{1}{2} \frac{|A|^2 |B|^2}{|F|-1} \\
&\geq \frac{1}{2} |A| |B|
\end{align*}
by our hypothesis $|A| |B| \leq \frac{1}{2} |F|$.  The claim \eqref{boost} then follows by the pigeonhole principle.
\end{proof}

We now recall the following sumset estimates (see e.g. \cite{ruzsa}, \cite{n}):

\begin{lemma}[Sumset estimates]\label{sumset}  Let $A, B$ be a non-empty finite subsets of an additive group such that $|A+B| \leq K \min(|A|, |B|)$.  Then we have 
$$ |A \pm A \pm A \ldots \pm A| \leq C K^C |A|$$
for any additive combination of $A$, where the constants $C$ depend on the length of this additive combination.
\end{lemma}

Next, we recall Gowers' quantitative formulation \cite{gowers} of the Balog-Szemeredi lemma \cite{balog}:

\begin{theorem}\label{bz}\cite{gowers}, \cite{bourgain:high-dim}
Let $A,B$ be finite subsets of an additive group with cardinality $|A|=|B|$, and let $G$ be a subset of $A \times B$ with cardinality
$$ |G| \geq |A| |B| / K$$
such that we have the bound
$$ | \{ a+b: (a,b) \in G \} | \leq K |A|.$$
Then there exists subsets $A'$, $B'$ of $A$ and $B$ respectively with $|A'| \geq c K^{-C} |A|$, $|B'| \geq c K^{-C} |B|$ such that
$$ |A'-B'| \leq C K^C |A|.$$
Indeed, we have the stronger statement that for every $a' \in A$ and $b' \in B$,
there are at least $cK^{-C} |A|^5$ solutions to the problem
$$ a'-b' = (a_1-b_1)-(a_2-b_2)+(a_3-b_3); a_1,a_2,a_3 \in A; b_1,b_2,b_3 \in B.$$
\end{theorem}

Note that all the above additive theorems have multiplicative analogues for multiplicative groups.  In particular there are multiplicative analogues on $F$ provided we eliminate the origin 0 from $F$ (though in our applications this single element is insignificant to our estimates).

We now recall a lemma from \cite{katz:falc} (see also \cite{katz:falc-expos},
\cite{bourgain:ring}):

\begin{lemma}\label{kt}\cite{katz:falc}  Let $A$ be a non-empty subset of $F$
such that
$$ |A+A|, |A \cdot A| \leq K |A|.$$
Then there is a subset $A'$ of $A$ with $|A'| \geq c K^{-C} |A|$ such that
$$ |A'\cdot A' - A' \cdot A'| \leq C K^C |A'|.$$
\end{lemma}

\begin{proof}  We outline the argument from \cite{katz:falc} or \cite{bourgain:ring}.  We shall use $X \lessapprox Y$ to
denote the estimate $X \leq C K^C Y$.  Without loss of generality we may assume that $|A| \gg 1$ is large, and that 
$0 \neq A$ (since removing $0$ from $A$ does not significantly affect any of the hypotheses).

We first observe from Theorem \ref{bz} that we can find subsets $C$, $D$ of $A$ with $|C|, |D| \approx |A|$ such that every element in $C-D$ has $\gtrapprox |A|^5$ representations of the form
$$ a_1 - a_2 + a_3 - a_4 + a_5 - a_6; \quad a_1, \ldots, a_6 \in A.$$
Multiplying this by an arbitrary element of $A \cdot A \cdot A / A \cdot A$, 
we see that every element of $(C-D) \cdot A \cdot A \cdot A / A \cdot A$
has $\gtrapprox |A|^5$ representations of the form
$$ b_1 - b_2 + b_3 - b_4 + b_5 - b_6; \quad b_1, \ldots, b_6 \in A \cdot A \cdot A \cdot A / A \cdot A.$$
However, by the multiplicative form of Lemma \ref{sumset}, the set $A \cdot A \cdot A \cdot A / A \cdot A$ has cardinality $\approx |A|$.  Thus by Fubini's theorem we have
\begin{equation}\label{cdaaaa}
 |(C-D) \cdot A \cdot A \cdot A / A \cdot A| \lessapprox |A|.
\end{equation}
Now we refine $C$ and $D$.  Since $|C|, |D| \approx |A|$ and $|A \cdot A| \approx |A|$, we have $|CD| \approx |C|, |D|$, and hence by the multiplicative form of Theorem \ref{bz}, we can find subsets $C'$, $D'$ of $C$, $D$ with $|C'|, |D'| \approx |A|$ such that every element in $C'D'$ has $\gtrapprox |A|^5$ representations in the form
$$ \frac{c_1 d_1 c_3 d_3}{c_2 d_2}; \quad c_1,c_2,c_3 \in C; \quad d_1,d_2,d_3 \in D.$$
Now let $c,c' \in C'$ and $d,d' \in D$ be arbitrary.  By the pigeonhole principle there thus exist $c_2 \in C$, $d_2 \in D$ such that we have $\gtrapprox |A|^3$ solutions to the problem
$$ cd = \frac{c_1 d_1 c_3 d_3}{c_2 d_2}; \quad c_1,c_3 \in C; \quad d_1,d_3 \in D.$$
We can rewrite this as
$$ cd - c'd' = x_1 - x_2 + x_3 + x_4$$
where
\begin{align*}
x_1 &= \frac{(c_1 - d') d_1 c_3 d_3}{c_2 d_2} \\
x_2 &= \frac{d'(c'-d_1) c_3 d_3}{c_2 d_2}\\
x_3 &= \frac{d'c'(c_3-d_2)d_3}{c_2d_2}\\
x_4 &= \frac{d'c'd_2(c_2-d_3)}{c_2d_2}.
\end{align*}
For fixed $c,d,c',d',c_2,d_2$, it is easy to see that the map from 
$(c_1,c_3,d_1,d_3)$ to $(x_1,x_2,x_3,x_4)$ is a bijection.  Since all the $x_j$ lie in $(C-D) \cdot A \cdot A \cdot A / A \cdot A$, we thus have $\gtrapprox |A|^3$ ways to represent $cd-c'd'$ in the form $x_1 - x_2 + x_3 - x_4$, where $x_1,x_2,x_3,x_4$ all lie in $(C-D) \cdot A \cdot A \cdot A / A \cdot A$.
By \eqref{cdaaaa} and Fubini's theorem we thus have
$$ |C'D' - C'D'| \lessapprox |A|.$$
In particular we have $|C'D'| \lessapprox |A| \lessapprox |C'|$, which by the multiplicative form of Lemma \ref{sumset} implies $|C'/D'| \approx |C'|$.  By considering the fibers of the quotient map $(x,y) \to x/y$ on $C' \times D'$ and using the pigeonhole principle, we thus see that there must be a non-zero field element $x$ such that $|C' \cap D'x| \approx |A|$.  If we then set $A' := C' \cap D'x$ we have $|A' A' - A' A'| \lessapprox |A|$ as desired.
\end{proof}

In the next section we bootstrap the $A \cdot A - A \cdot A$ type bound in Lemma \ref{kt} to other polynomial expressions of $A$.

\section{Iterated sum and product set estimates}

We now prove the following lemma, which is in the spirit of Lemma \ref{sumset}:

\begin{lemma}\label{iterated}  Let $A$ be a non-empty subset of a finite field $F$, and suppose that we have the bound
$$ |A.A - A.A| \leq K|A|$$
for some $K \geq 1$.  We adopt the normalization that $1 \in A$.
Then for any polynomial $P$ of several variables and integer coefficients, we have
$$ |P(A,A,\ldots,A)| \leq C K^C |A|$$
where the constants $C$ depend of course on $P$.
\end{lemma}

\begin{proof}  We need some notation.  We say that a set $A$ is \emph{essentially contained} in $B$, and write $A \Subset B$, if we have $A \subseteq X + B$ for some set $X$ of cardinality $|X| \leq C K^C$.

We have the following simple lemma of Ruzsa \cite{ruzsa-group}:

\begin{lemma}\label{ruzsa}  Let $A$ and $B$ be subsets of $F$ such that $|A+B| \leq CK^C |A|$ or $|A-B| \leq CK^C |A|$.  Then $B \Subset A-A$.
\end{lemma}

\begin{proof}  By symmetry we may assume that $|A+B| \leq CK^C |A|$.  Let $X$ be a maximal subset of $B$ with the property that the sets $\{ x+A: x \in X\}$ are all disjoint.  Since the sets $x+A$ all have cardinality $|A|$ and are all contained in $A+B$, we see from disjointness that $|X| |A| \leq |A+B|$, and hence $|X| \leq CK^C$.  Since the set $X$ is maximal, we see that for every $b \in B$, the set $b+A$ must intersect $x+A$ for some $x \in X$.  Thus $b \in x+A-A$, and hence $B \subseteq X + A-A$ as desired.
\end{proof}

Call an element $x \in F$ \emph{good} if we have $x \cdot A \Subset A-A$.

\begin{proposition}\label{it}  The following three statements are true.
\begin{itemize}
\item Every element of $A$ is good.
\item If $x$ and $y$ are good, then $x+y$ and $x-y$ is good.
\item If $x$ and $y$ are good, then $xy$ is good.
\end{itemize}
(Of course, the implicit constants in ``good'' vary at each occurence).
\end{proposition}

\begin{proof}  Let us first show that every element of $A$ is good.  Since $1 \in A$, we have
$$ |A.A - A| \leq |A.A-A.A| \leq K|A|$$
and hence by Lemma \ref{ruzsa}
\begin{equation}\label{aaa}
A.A \Subset A-A
\end{equation}
which implies in particular that every element of $x$ is good.

Now suppose that $x$ and $y$ are good, thus $x \cdot A \Subset A-A$ and $y \cdot A \Subset A-A$.  Then
$$ (x+y) \cdot A \subseteq x \cdot A + y \cdot A \Subset A-A + A-A.$$
On the other hand, since $|A-A| \leq |A.A-A.A| \leq K|A|$, we have from sumset estimates (Lemma \ref{sumset}) that 
$$ |A-A+A-A+A| \leq  CK^C |A|$$
and hence by Lemma \ref{ruzsa}
\begin{equation}\label{a+a}
A-A+A-A \Subset A-A.
\end{equation}
Thus by transitivity of $\Subset$ we have $(x+y) \cdot A \Subset A-A$ and hence $x+y$ is good.  A similar argument shows that $x-y$ is good.

Now we need to show that $xy$ is good.  Since $x \cdot A \Subset A-A$ we have
$$ xy \cdot A \Subset y \cdot A - y \cdot A.$$
But since $y \cdot A \Subset A-A$, we have
$$ xy \cdot A \Subset A-A - A + A.$$
By \eqref{a+a} we conclude that $xy$ is good.
\end{proof}

By iterating this proposition we see that for any integer polynomial $P$, every element of $P(A,\ldots,A)$ is good\footnote{An alternate way to proceed at this point is to show that the number of good points is at most $\lessapprox N$; indeed, it is easy to show that any good point is contained inside $(A-A+A-A)/(A-A)$ if $N$ is sufficiently large, where we exclude 0 from the denominator $A-A$ of course.  We omit the details.}.

Write $A^2 := A \cdot A$, $A^3 := A \cdot A \cdot A$, etc.  We now claim inductively that $A^k \Subset A-A$ for all $k=0,1,2,3,\ldots$.  The case $k=0,1$ are trivial, and $k=2$ has already been covered by
\eqref{aaa}.  Now suppose inductively that $k>2$, and that we have already proven that 
$A^{k-1} \Subset A-A$.  Thus 
$$ A^{k-1} \subseteq X + A - A$$
for some set $X$ of cardinality $|X| \leq C K^C$.  Clearly we may restrict $X$ to the set $A^{k-1} - (A-A)$.  In particular, every element of $X$ is good.  We now multiply by $A$ to obtain
$$ A^k \subseteq X \cdot A + A \cdot A - A \cdot A.$$
Since every element of $X$ is good, and $|X| \leq CK^C$, we see that $X \cdot A \Subset A-A$.  By \eqref{aaa} we thus have
$$ A^k \Subset A-A+A-A-(A-A).$$
But by arguing as in the proof of \eqref{a+a} we have
$$ A-A+A-A-(A-A) \Subset A-A,$$
and thus we can close the induction.

Since $A^k \Subset A-A$ for every $k$, and $A-A \pm (A-A) \Subset A-A$ by \eqref{a+a}, we thus 
see that every integer combination of $A^k$ is essentially contained in $A-A$.  In particular 
$P(A,\ldots,A) \Subset A-A$ for every integer polynomial $A$, and the claim follows.
\end{proof}

\section{Proof of the sum-product estimate}\label{sp-sec}

We now have all the machinery needed to prove Theorem \ref{sum-product}.
We basically follow the Edgar-Miller approach, see \cite{edgar}.  We write $F$ for $\Z/q\Z$, and let $F^* := F - \{0\}$ be the invertible elements of $F$.  Let $\delta > 0$, and let $A$ be a subset of $F$ such that $|F|^\delta < |A| < |F|^{1-\delta}$.

Let $0 < \eps \ll 1$ be a small number depending on $\delta$ to be chosen later.  In this section we use $X \lesssim Y$ to denote the estimate $X \leq C(\delta,\eps) Y$ for some $C(\delta,\eps) > 0$.  Suppose for contradiction that
$$ |A+A|, |A \cdot A| \lesssim |A|^{1+\eps};$$
Then by Lemma \ref{kt}, and passing to a refinement of $A$ if necessary, we may assume that
$$ |A \cdot A - A \cdot A|  \lesssim |A|^{1+C \eps}.$$
We may normalize $1 \in A$.  By Lemma \ref{iterated} we thus have
\begin{equation}\label{polynomial}
 |P(A, \ldots, A)| \lesssim |A|^{1+C\eps}
\end{equation}
for any polynomial $P$ with integer coefficients, where the constants $C$
depend of course on $P$.

Our first objective is to obtain a linear surjection from $A^k$ to $F$ for sufficiently large $k$:

\begin{lemma}\label{ak}
There exists a positive integer $k \sim 1/\delta$, and invertible
field elements $\xi_1, \ldots, \xi_k \in F^*$, such that
$$
F =A\xi_1+ \cdots+ A\xi_k.
$$
In other words, we have a linear surjection from $A^k$ to $F$.
\end{lemma}

\begin{proof}
Iterating Lemma \ref{mixed} about $O(1/\delta)$ times, we obtain $\xi_1, \ldots, \xi_k \in F^*$ such that
$$
|A\xi_1+\cdots+A\xi_k|\geq\frac{|F|}{10}.
$$
The lemma then obtains after $O(1)$ applications of the Cauchy-Davenport inequality \eqref{cauchy-davenport},
increasing $k$ as necessary.
\end{proof}

Next, we reduce the rank $k$ of this surjection, at the cost of replacing $A$
by a polynomial expression of $A$.

\begin{lemma}\label{iterate}  Let $B$ be a non-empty subset of $F$, and suppose $k > 1$ is such that there is a linear surjection from $B^k$ to $F$.  Then there is a linear surjection from $\tilde B^{k-1}$ to $F$, where $\tilde B := B \cdot(B-B) + B \cdot (B-B)$.
\end{lemma}

\begin{proof}
By hypothesis, we have a surjection
$$B^k\rightarrow F: (a_1, \ldots, a_k)\mapsto \sum_{j\leq k} a_j\xi_j$$
for some $\xi_1, \ldots, \xi_k \in F$.
Our map cannot be one-to-one, since otherwise
$$
|B|^k = |F| \text { (contradicting primarily of $|F|$)}.
$$
Thus there are $(b_1, \ldots, b_k) \not= (b_1', \ldots b_k') \in B^k$ with
\begin{equation}\label{21}
(b_1-b_1')\xi_1+\cdots+ (b_k-b_k')\xi_k=0.
\end{equation}
Let $b_k\not= b_k'$.  By the surjection property
$$ F = B \xi_1 + \ldots + B \xi_k;$$
since $F$ is a field, we thus have
$$
F=B\xi_1(b_k-b_k')+\cdots+ B\xi_k(b_k-b_k')
$$
and substituting $(b_k-b_k')\xi_k$ from \eqref{21}
\begin{align*}
F &= B\xi_1(b_k-b_k')+\cdots+ B\xi_{k-1}(b_k-b_{k}')-B(b_1-b_1')\xi_1-\cdots- B(b_{k-1}-b_{k-1}') \xi_{k-1}\\
&\subset \tilde B\xi_1 +\cdots+ \tilde B\xi_{k-1}
\end{align*}
and the claim follows.
\end{proof}

Starting with Lemma \ref{ak} and then iterating Lemma \ref{iterate} $k$ times, we eventually get a linear 
surjection from a polynomial expression $P(A,\ldots,A)$ of $A$ to $F$, and thus
$$ |P(A,\ldots,A)| \geq |F|.$$
But this contradicts \eqref{polynomial}, if $\eps$ is sufficiently small
depending on $\delta$.  This contradiction proves Theorem \ref{sum-product}.
\endprf

{\bf Remark.}  Suppose the finite field $F$ did not have prime order.  Then the analogue of Theorem \ref{sum-product}
fails, since one can take $A$ to be a subfield $G$ of $F$, or a large subset of such a subfield $G$.  It turns
out that one can adapt the above argument to show that these are in fact the only ways in which Theorem \ref{sum-product}
can fail (up to dilations, of course):

\begin{theorem}\label{general}  Let $A$ be a subset of a finite field $F$ such that $|A| > |F|^\delta$ for some 
$0 < \delta < 1$, and suppose that $|A+A|, |A \cdot A| \leq K |A|$ for some $K \gg 1$.  Then there exists a subfield $G$ of $F$ of cardinality $|G| \leq K^{C(\delta)} |A|$, a non-zero field element $\xi \in F - \{0\}$, and a set $X \subseteq F$ of cardinality $|X| \leq K^{C(\delta)}$ such that $A \subseteq \xi G \cup X$.
\end{theorem}

It is interesting to compare the above theorem to Freiman's theorem (\cite{freiman}, \cite{ruzsa}, \cite{chang}) which does not assume
control on $|A \cdot A|$ but has a dependence on constants which is significantly worse than polynomial.  It seems
possible that the constant $C(\delta)$ can be made independent of $\delta$, but we do not know how to do so.

\begin{proof}  (Sketch) Of course, we may assume that $|F| \geq K^{C(\delta)}$ for some large $C(\delta)$.  We repeat the argument used to prove Theorem \ref{sum-product}.  This argument allows us to find a refinement $A'$ of $A$ with $|A'| \geq K^{-C} A$ such that $|A' \cdot A' - A' \cdot A'| \leq K^C |A|$.  By dilating $A$ and $A'$ if necessary we may assume as before that $1 \in A'$ (as we shall see, this normalization allows us to take $\xi = 1$ in the conclusion of this theorem).  By Lemma \ref{iterated} we thus have $|P(A',\ldots,A')| \leq K^C |A'|$ for all integer polynomials $P$, 
with the 
constant $C$ depending on $P$ of course.  We may assume $0 \in A'$ since adding 0 to $A'$ and $A$ do not significantly affect the above
polynomial bounds.  

We now claim that $A'$ is contained in some subfield $G$ of $F$ of cardinality $|G| \leq K^{C(\delta)} |A|$.  The argument in Lemma \ref{ak} still gives a 
surjection from $(A')^k$ to $F$ for some $k \sim 1/\delta$.  We then attempt to use Lemma \ref{iterate} to drop the rank of this surjection down to 1.
If we can reduce the rank all the way to one, then we have by arguing as before that $|F| \leq K^{C(k)} |A|$, so the claim follows by
setting $G := F$.  The only time we run into difficulty in this iteration is if we discover a linear 
surjection from some $\tilde A^{k'}$ to $F$ with $k' > 1$ which is also injective, where $\tilde A$ is some 
polynomial expression of $P(A',\ldots,A')$.  An inspection of the proof of Lemma \ref{iterate}, combined with the normalizations $0,1 \in A'$, reveals that $\tilde A$ must contain $A'$.  If we have 
$|\tilde A + \tilde A| > |\tilde A|$, then the linear map from $(\tilde A + \tilde A)^{k'} \to F$
is surjective but not injective, which allows us to continue the iteration of Lemma \ref{iterate}.  Similarly
if $|\tilde A \cdot \tilde A| > |\tilde A|$.  Thus the only remaining case is when $|\tilde A| = |\tilde A+ \tilde A| 
= |\tilde A \cdot \tilde A|$.  But this, combined with the fact that $0,1 \in \tilde A$, implies that 
$\tilde A = \tilde A + \tilde A = \tilde A \cdot \tilde A$,
and hence that $\tilde A$ is a subfield of $F$.  Since $|\tilde A| \leq K^{C(k)} |A'|$, the claim follows.

This shows that $A'$ is a subset of $G$.  Since $|A+A'| \leq K|A|$, we see from Lemma \ref{ruzsa} that $A \Subset A' - A'$, and hence $A \Subset G$.  Thus there exists a set $Y$ of cardinality $|Y| \leq K^{C(\delta)}$ such that $A \subseteq G + Y$.

Let $y \in Y - G$.  To finish the proof (with $\xi = 1$) it will suffice to show that $|A \cap (G + y)| \leq K^{C(\delta)}$ for all such $y$.  But observe that for any two distinct $x, x' \in G+y$, the sets $xG$ and $x'G$ do not intersect except at the origin (for if $xg = x'g'$, then $g \neq g'$, and hence $x = (x'-x) \frac{g'}{g-g'} \in G$, contradicting the hypotheses that $x \in G+y$ and $y \not \in G$).  In particular, the sets $x(A'-\{0\})$ and $x'(A'-\{0\})$ are disjoint.  Thus
$$ K|A| \geq |A \cdot A| \geq  |A \cap (G + y)| |A' - \{0\}| \geq |A \cap (G+y)| K^{-C} |A|$$
and the claim follows.

\end{proof}

\section{Some basic combinatorics}\label{notation-sec}

In later sections we shall use the sum-product estimate in Theorem \ref{sum-product} to various combinatorial
problems in finite geometries.  In doing so we will repeatedly use a number of basic combinatorial tools,
which we collect here for reference.

We shall frequently use the following elementary observation:  If $B$ is a finite set, and $\mu: B \to \R^+$ is a function such that
$$ \sum_{b \in B} \mu(b) \geq X,$$
then we have

$$ \sum_{b \in B: \mu(b) \geq X/2|B|} \mu(b) \geq X/2.$$
We refer to this as a ``popularity'' argument, since we are restricting $B$ to the values $b$ which are ``popular'' in the sense that $\mu$ is large.  

We shall frequently use the following version of the Cauchy-Schwarz inequality.

\begin{lemma}\label{cz}  Let $A$, $B$ be finite sets, 
and let $\sim$ be a relation connecting pairs $(a,b) \in A \times B$ such that
$$ |\{ (a,b) \in A \times B: a \sim b \}| \gtrsim X$$
for some $X \gg |B|$.  Then
$$ |\{ (a,a',b) \in A \times A \times B: a \neq a'; a, a' \sim b \}| \gtrsim \frac{X^2}{|B|}.$$
\end{lemma}

\begin{proof}
Define for each $b \in B$, define $\mu(b) := | \{ a \in A: a \sim b \} |$.  Then by hypothesis we have
$$ \sum_{b \in B} \mu(b) \gtrsim X.$$
In particular, by the popularity argument we have
$$ \sum_{b \in B: \mu(b) \gtrsim X/|B|} \mu(b) \gtrsim X.$$
By hypothesis, we have $X/|B| \gg 1$.  From this and the previous, we obtain
$$ \sum_{b \in B: \mu(b) \gtrsim X/|B|} \mu(b) (\mu(b)-1) \gtrsim X (X/|B|)$$
and the claim follows.
\end{proof}

A typical application of the above Lemma is the standard incidence bound on lines in a plane $F^2$, where $F$ is a finite field.  

\begin{corollary}\label{easy-cor}  Let $F^2$ be a finite plane.  For an arbitrarily collection $P \subseteq F^2$ of points and $L$ of lines in $F^2$, we have
\begin{equation}\label{easy-incidence}
|\{ (p,l) \in P \times L: p \in l \}| \lesssim |P|^{1/2} |L| + |P| 
\end{equation}
\end{corollary}

\begin{proof}
We may of course assume that the left-hand side of \eqref{easy-incidence} is $\gg |P|$, since the claim is trivial otherwise.  From Lemma \ref{cz} we have
$$ |\{ (p, l, l') \in P \times L \times L: p \in l \cap l'; l \neq l' \}| \gtrsim
|P|^{-1} |\{ (p,l) \in P \times L: p \in l \}|^2.$$
On the other hand, $|l \cap l'|$ has cardinality $O(1)$ if $l \neq l'$, thus
$$ |\{ (p, l, l') \in P \times L \times L: p \in l \cap l'; l \neq l' \}| \lesssim 
|L|^2.$$
Combining the two estimates we obtain the result.
\end{proof}

\section{A Szemer\'edi-Trotter type theorem in finite fields}\label{sz-sec}

We now use the one-dimensional sum-product estimate to obtain a key two-dimensional estimate, namely an incidence bound of Szemer\'edi-Trotter type. 

Let $F$ be a finite field, and consider the projective finite plane $PF^3$, which is the set $F^3 - \{(0,0,0)\}$
quotiented by dilations.  We embed the ordinary plane $F^2$ into $PF^3$ by identifying $(x,y)$ with the
equivalence class of $(x,y,1)$; $PF^3$ is thus $F^2$ union the line at infinity.
Let $1 \leq N \leq |F|^2$ be an integer, and let $P$ be a collection of points and $L$ be a collection of lines in $F^2$.  We consider the problem of obtaining an upper bound on the number of incidences
$$ |\{ (p,l) \in P \times L: p \in l \}|.$$
From Corollary \ref{easy-cor} and the duality between points and lines in two dimensions we have the easy bounds 
\begin{equation}\label{incidence}
 |\{ (p,l) \in P \times L: p \in l \}| \leq \min(|P| |L|^{1/2} + |L|, |L| |P|^{1/2} + |P|),
\end{equation}
see e.g. \cite{bollobas}.  In a sense, this is sharp:
if we set $N = |F|^2$, and let $P$ be all the points in $F^2 \subset PF^3$ and $L$ be most of the lines in $F^2$, then we have roughly $|F|^3 \sim N^{3/2}$ incidences.  More generally if $G$ is any subfield of $F$ then one can construct a similar example with $N = |G|^2$,
$P$ being all the points in $G^2$, and $L$ being the lines with slope and intercept in $G$.  

Recently Elekes \cite{elekes} observed that there is a connection between this incidence problem and the sum-product problem:

\begin{lemma}\cite{elekes} Let $A$ be a subset of $F$.  Then there is a collection of points $P$ and lines $L$ with $|P| = |A+A| |A \cdot A|$ and $|L| = |A|^2$ which has at least $|A|^3$ incidences.
\end{lemma}

\begin{proof} Take $P = (A+A) \times (A\cdot A)$, and let $L$ be the set of all lines of the form $l(a,b) := \{ (x,y): y = b(x-a) \}$ 
where $a,b$ are any two elements of $A$.  The claim follows since $(a+c,bc) \in P$ is incident to $l(a,b)$ whenever $a,b,c \in A$.
\end{proof}

Thus any improvement to the trivial bound of $O(N^{3/2})$ on the incidence problem should imply a sum-product estimates.  Conversely,
we can use the sum-product estimate (Theorem \ref{sum-product}) to obtain a non-trivial incidence bound:

\begin{theorem}\label{sz} Let $F$ be the finite field $F := \Z/q\Z$ for some prime $q$, and let $P$ and $L$ be points and lines in $PF^3$ with cardinality $|P|,|L| \leq N=|F|^\alpha$ for some $0 < \alpha < 2$.  Then we have
$$  |\{ (p,l) \in P \times L: p \in l \}| \leq C N^{3/2 - \eps} $$
for some $\eps = \eps(\alpha) > 0$ depending only on the exponent $\alpha$.
\end{theorem}

{\bf Remark.}  The corresponding statement for $N$ points and $N$ lines in the Euclidean plane $\R^2$ (or $P\R^3$) is due to Szem\'eredi and Trotter \cite{Szemeredi:trotter}, with $\eps := 1/6$.  This bound is sharp.  It may be that one could similarly take $\eps = 1/6$ in the finite field case when $\alpha$ is sufficiently small, but we do not know how to do so; certainly the argument in \cite{Szemeredi:trotter} relies crucially on the ordering properties of $\R$ and so does not carry over to finite fields.

\begin{proof}  We may assume that $N \gg 1$ is large.  By adding dummy points and lines we may assume that $|P|=|L|=N$.

Fix $N = |F|^\alpha$, and let $0 < \eps \ll 1$, be chosen later.  Suppose for contradiction that we can find points 
$P$ and lines $L$ with $|P|=|L|=N$ such that
$$  |\{ (p,l) \in P \times L: p \in l \}| \gtrsim N^{3/2 - \eps};$$
we shall use the sum-product estimates to obtain a contradiction if $\eps$ is sufficiently small.  Our arguments follow those
in \cite{katz:falc}, \cite{katz:falc-expos}.

We first use the popularity argument to control how many points are incident to a line and vice versa.  For each $p \in P$, define the multiplicity $\mu(p)$ at $p$ by
$$ \mu(p) := |\{ l \in L: p \in l \}|.$$
Then by hypothesis
$$ \sum_{p \in P} \mu(p) \gtrsim N^{3/2-\eps}$$
and hence by the popularity argument and the hypothesis $|P|=N$

$$ \sum_{p \in P: \mu(p) \gtrsim N^{1/2-\eps}} \mu(p) \gtrsim N^{3/2-\eps}.$$
On the other hand, we observe that
\begin{align*}
\sum_{p \in P: \mu(p) \gg N^{1/2+\eps}} \mu(p)^2
&\ll N^{-1/2-\eps}
\sum_{p \in P} \mu(p) (\mu(p)-1)\\
&\lesssim
N^{-1/2-\eps} \sum_{p \in P} |\{(l, l') \in L \times L: p \in l,l'; l \neq l' \}| \\
& = N^{-1/2-\eps} \sum_{l,l' \in L: l \neq l'} |\{p \in P: p \in l,l' \}|\\
&\leq N^{-1/2-\eps} \sum_{l,l' \in L: l \neq l'} 1 \\
&\leq N^{1/2+\eps}.
\end{align*}
Thus if we set $P' \subseteq P$  to be the set of all points $p$ in $P$ such that
$$ N^{1/2-\eps} \lesssim \mu(p) \lesssim N^{1/2+\eps}$$
then we have
$$ \sum_{p \in P'} \mu(p) \gtrsim N^{3/2-\eps}.$$

For each $l \in L$, define the multiplicity $\lambda(l)$ by
$$ \lambda(l) := |\{ p \in P': p \in l \}|,$$
then we can rewrite the previous as
$$ \sum_{l \in L} \lambda(l) \gtrsim N^{3/2-\eps}.$$
By the popularity argument we thus have
$$ \sum_{l \in L: \lambda(l) \gtrsim N^{1/2-\eps}} \lambda(l) \gtrsim N^{3/2-\eps}.$$
On the other hand, we have
\begin{align*}
\sum_{l \in L: \lambda(l) \gg N^{1/2+\eps}} \lambda(l) 
&\lesssim N^{-1/2-\eps}
\sum_{l \in L} \lambda(l) (\lambda(l)-1) \\
&\lesssim N^{-1/2-\eps}
\sum_{l \in L} |\{ (p,p') \in P' \times P': p,p' \in l; p \neq p' \}| \\
&=
N^{-1/2-\eps}
\sum_{p,p' \in P': p \neq p'} |\{ l \in L: p,p' \in l \}| \\
&\leq
N^{-1/2-\eps}
\sum_{p,p' \in P': p \neq p'} 1.
\end{align*}
Thus if we set $L' \subset L$ to be the set of all lines $l$ in $L$ such that
$$ N^{1/2-\eps} \lesssim \lambda(p) \lesssim N^{1/2+\eps}$$
then we have
$$ \sum_{l \in L'} \lambda(l) \gtrsim N^{3/2-\eps}.$$
For each $p \in P'$, let $\mu'(p)$ denote the multiplicity
$$ \mu'(p) := | \{ l \in L': p \in l \}|;$$
clearly $\mu'(p) \leq \mu(p)$.  We can then rewrite the previous estimate as
$$ \sum_{p \in P'} \mu'(p) \gtrsim N^{3/2-\eps}.$$
Thus by the popularity argument, if we set $P'' \subseteq P'$ to be the set of all points $p$ in $P'$ such that 
$$ \mu'(p) \gtrsim N^{1/2-\eps}$$
then we have
$$ \sum_{p \in P''} \mu'(p) \gtrsim N^{3/2-\eps}.$$
or equivalently 
$$  |\{ (p,l) \in P'' \times L': p \in l \}| \gtrsim C_0 N^{3/2 - \eps}.$$
Since $|L'| \leq N$, we have in particular that
\begin{equation}\label{pp}
|P''| \gtrsim N^{1/2-\eps}.
\end{equation}

The next step is to capture a large portion of the popular point set $P'$ inside a Cartesian product $A \times B$,
possibly after a projective transformation.  The key observation is that such a product arises, modulo projective
transformations, whenever one intersects two ``bushes'' of lines.

Let $p_0$ be any point in $P''$.  Then by construction there are $\gtrsim N^{1/2-\eps}$ lines $l$ in $L'$ containing $p_0$.
Each of these lines $l$ contains $\gtrsim N^{1/2-\eps}$ points $p$ in $P'$;
of course, all but one of these are distinct from $p_0$.  Thus we have
$$ |\{ (p,l) \in P' \times L': p,p_0 \in l, p \neq p_0 \}| \gtrsim N^{1-2\eps}.$$
Let us define a relation $\sim$ on $P$ by defining $p \sim p'$ if $p \neq p'$ and there is a line in $L'$ containing both $p$ and $p'$.  Since two distinct points determine at most one line, we thus have
$$ | \{ p \in P': p \sim p_0 \}| \gtrsim N^{1-2\eps} \hbox{ for all } p_0 \in P''.$$
Summing this over all $p_0$ in $P''$, we obtain
$$ |\{ (p_0,p) \in P'' \times P': p \sim p' \}| \gtrsim N^{1-2\eps} |P''|.$$

Since $|P'| \leq N$, we thus see by Lemma \ref{cz} that
$$ |\{ (p_0,p_1,p) \in P'' \times P'' \times P': p \sim p_0,p_1; p_0 \neq p_1 \}| \gtrsim N^{1-C\eps} |P''|^2.$$
By the pigeonhole principle, there thus exist distinct points $p_0,p_1 \in P''$ such that
\begin{equation}\label{proj}
 |\{ p \in P': p \sim p_0,p_1\}| \gtrsim N^{1-C\eps}.
\end{equation}
Fix these $p_0,p_1$.  By applying a projective linear transformation (which maps lines to lines and preserves
incidence) we may assume that $p_0$, $p_1$ are both on the line at infinity.  Indeed, we may assume that 
$p_0 = [(1,0,0)]$ and $p_1 = [(0,1,0)]$, where $[(x,y,z)]$ is the equivalence class of $(x,y,z)$ in $PF^3$.

We first eliminate those points $p$ in \eqref{proj} on the line at infinity.  Such points can only occur if the 
line at infinity is in $L'$.  But then that line contains at most $O(N^{1/2+\eps})$ points in $P'$, by the definition of $L'$.  Thus 
if $\eps$ is sufficiently small we have
$$  |\{ p \in P' \cap F^2: p \sim [(1,0,0)], [(0,1,0)]\}| \gtrsim N^{1-C\eps}.$$
Consider the lines in $L'$ which pass through $[(1,0,0)]$.  In the plane $F^2$, these lines be horizontal, i.e.
they are of the form $\{ (x,y) \in F^2: y = b \}$ for some $b \in F$.  Let $B \subseteq F$ denote the set of all 
such $b$.  Since each line contains at least $c N^{1/2-\eps}$ points in $P'$, and $|P'| \leq N$, we know that
$|B| \lesssim N^{1/2+\eps}$.  Similarly the lines in $L'$ which pass through $[(0,1,0)]$ must in $F^2$
be vertical lines of the form $\{ (x,y) \in F^2: x = a \}$ for $a \in A$, where $|A| \leq C N^{1/2+\eps}$.
We thus have
\begin{equation}\label{p-large}
 |P' \cap (A \times B)| \gtrsim N^{1-C\eps}
\end{equation}
and
\begin{equation}\label{ab-small} |A|, |B| \leq C N^{1/2 + \eps}
\end{equation}

Now that we have placed $P'$ in a Cartesian grid, the next step is to exploit the form $y=mx+b$ of lines
in $F^2$ to obtain some additive and multiplicative information on $A$ and $B$.

Define $P_0 := P' \cap (A \times B)$.  By definition of $P'$ we have

$$ |\{ l \in L: p \in l \}| \gtrsim N^{1/2-\eps} \hbox{ for all } p \in P_0;$$
summing over $P_0$ using \eqref{p-large} and rearranging, we obtain
$$ | \{ (p,l) \in P_0 \times L: p \in l \}| \gtrsim N^{3/2 - C\eps}.$$
Let $L_0$ be those lines in $L$ which are not horizontal.  Since horizontal lines can contribute at most $|P_0| \leq N$ incidences to the above expression,
we have (if $\eps$ is sufficiently large)
$$ | \{ (p,l) \in P_0 \times L_0: p \in l \}| \gtrsim N^{3/2 - C\eps}.$$
By the popularity argument, if we let $L_1$ denote those lines in $L_0$ such that
$$ | \{ p \in P_0: p \in l \} | \gtrsim N^{1/2 - C\eps}$$
we thus have
$$ | \{ (p,l) \in P_0 \times L_1: p \in l \} | \gtrsim N^{3/2 - C\eps}$$
if the implicit constants are chosen appropriately.

Define a relation $\sim$ between $B$ and $L_1$ by defining $b \sim l$ if there is a point
$p$ in the row $P_0 \cap (A \times \{b\})$ such that $p \in l$.  Note that such a point $p$ is unique since $l$
is not horizontal, and thus
$$ | \{ (b, l) \in B \times L_1: b \sim l \}| \gtrsim N^{3/2 - C\eps}.$$
By Lemma \ref{cz}, we thus have
$$ |\{(b,b',l) \in B \times B \times L_1: b,b' \sim l \}| \gtrsim N^{2 - C\eps}.$$
By \eqref{ab-small} and the pigeonhole principle, we thus conclude that there exists distinct heights $b,b' \in B$ such that
$$ |\{ l \in L_1: b,b' \sim l \}| \gtrsim N^{1-C\eps}.$$
Fix this $b,b'$.  By an affine transformation of the vertical variable 
(which does not affect the line at infinity) we may assume that
$b=0$ and $b'=1$.  Since each line $l \in L_1$ contains $\gtrsim N^{1/2-C\eps}$ points $(x,t)$ in $P_0$, and hence in 
$A \times B$, and most of these have $t \neq 0,1$ since $l$ is not horizontal, we have
$$ | \{ (x,t,l) \in A \times B \times L_1: 0, 1 \sim l; (x,t) \in l; t \neq 0,1 \} | \gtrsim N^{3/2 - C\eps}.$$
By definition of the relation $a \sim l$, we thus have
$$ | \{ (x,t,l,x_0,x_1) \in A \times B \times L_1 \times A \times A: (x_0,0), (x,t), (x_1,1) \in l; t \neq 0,1 \} |
\gtrsim N^{3/2 - C\eps}.$$

Since the three points $(x_0,0)$, $(x,t)$, $(x_1,1)$ determine $l$, and
$$ x = x_0 + (x_1-x_0) t,$$
we thus have
\begin{equation}\label{abb-high-0}
 | \{ (t,x_0,x_1) \in B \times A \times A: (1-t)x_0 + t x_1 \in A; t \neq 0,1 \}| \gtrsim N^{3/2 - C\eps}.
\end{equation}

Note that this is somewhat similar to saying that $(1-B).A + B.A \subseteq A$, so we are getting close
to being able to apply our sum-product estimate.  But first we must perform some Balog-Szemer\'edi
type refinements.

Let $A' \subseteq A$ denote those $x_1$ in $A$ for which
$$ | \{ (t,x_0) \in B \times A \times A: (1-t)x_0 + t x_1 \in A; t \neq 0,1 \}| \gtrsim N^{1 - C\eps}.$$
From \eqref{ab-small}, \eqref{abb-high-0} and the popularity argument we have 
\begin{equation}\label{abb-high}
 | \{ (t,x_0,x_1) \in B \times A \times A': (1-t)x_0 + t x_1 \in A; t \neq 0,1 \}| \gtrsim N^{3/2 - C\eps}
\end{equation}
if the implicit constants are chosen correctly.

In particular, from \eqref{ab-small} again we have
\begin{equation}\label{ap-card}
 |A'| \gtrsim N^{3/2 - C \eps} / |A||B| \gtrsim N^{1/2 - C \eps}.
\end{equation}
Also, by \eqref{abb-high}, the pigeonhole principle and \eqref{ab-small} we may find $t_0 \in B$ such that $t_0 \neq 0,1$ and
$$ | \{ (x_0,x_1) \in A \times A': (1-t_0)x_0 + t_0 x_1 \in A \}| \gtrsim N^{1/2 - C \eps}$$
By \eqref{ab-small} we have
$$ | \{ (x_0,x_1) \in A \times A': (1-t_0)x_0 + t_0 x_1 \in A \}| \gtrsim N^{-C\eps} |A| |A'|.$$
By \eqref{ap-card}, \eqref{ab-small} and Theorem \ref{bz} 
applied to the sets $(1-t_0)A$ and $t_0 A'$, we thus have a subsets $(1-t_0) \tilde A$ of $(1-t_0) A$ and
$t_0 A''$ of $t_0 A'$ with cardinalities at least $\gtrsim N^{1/2 - C\eps}$ such that
$$ |(1-t_0)\tilde A + t_0 A'| \lesssim N^{1/2 + C \eps}.$$
By \eqref{ap-card}, \eqref{ab-small} and sumset estimates, this implies in particular that
$$ |t_0 A' + t_0 A'| \lesssim N^{1/2 + C\eps}$$
and hence
\begin{equation}\label{ap-sum}
 |A' + A'| \lesssim N^{1/2 + C \eps}.
\end{equation}
Now we return to \eqref{abb-high}.  From \eqref{ab-small} and the pigeonhole principle we may find an $x_0 \in A$ such
that
$$
 | \{ (t,x_1) \in B \times A': (1-t)x_0 + t x_1 \in A; t \neq 0,1 \}| \gtrsim N^{1 - C \eps}.$$
By a translation in the horizontal variables $x_0$, $x_1$, $A$, $A'$ we may assume that $x_0 = 0$.  Thus
$$ | \{ (t,x_1) \in (B \backslash \{0\}) \times (A' \backslash \{0\}): t x_1 \in A\}| \gtrsim N^{1 - C \eps},$$
since the contribution of 0 is easily controlled by \eqref{ab-small}.  By \eqref{ab-small} and the multiplicative form of Theorem \ref{bz}, we can thus find a subset $A''$ of $A' \backslash \{0\}$ with $|A''| \gtrsim N^{1/2-C\eps}$ and
$$ |A'' \cdot A''| \lesssim N^{1/2 + C \eps}.$$
On the other hand, from \eqref{ap-sum} we have
$$ |A'' + A''| \lesssim N^{1/2 + C \eps}.$$
But this gives a contradiction to the sum product estimate (Theorem \ref{sum-product}) if $\eps$ is sufficiently small.
\end{proof}

{\bf Remark.}  One can extend this result to more general finite fields $F$ using Theorem \ref{general} as a substitute
for Theorem \ref{sum-product}.  Informally,
the result is as follows: for general finite fields, one has the same conclusions as Theorem \ref{sz} except
when $P$ has large intersection with a projective transformation of a Cartesian product $G \times G$ for some
subfield $G$ of $F$, and when $L$ has large intersection with the associated collection of lines.  We omit the details.  We also remark that similar generalizations can be made for the problems stated in the next few sections, but the generalizations become quite cumbersome to state and prove and we shall not do so here.

\section{Applications to the distance set problem}\label{distance-sec}

We now work in the finite field plane $F^2$.  Given any two points $(x_1,y_1)$, $(x_2,y_2)$, we define the distance 
$d((x_1,y_1),(x_2,y_2)) \in F$ by
$$ d((x_1,y_1), (x_2,y_2)) = (x_1-x_2)^2 + (y_1 - y_2)^2$$
(we omit the square root to avoid some distracting technicalities).  Given any collection $P$ of points in $F^2$, we define the distance set $\Delta(P) \subseteq F$ by
$$ \Delta(P) := \{ d(p,p'): p,p' \in P \}.$$
The \emph{Erd\"os distance problem} is to obtain the best possible lower bound for $|\Delta(P)|$ in terms of $|P|$.  If $-1$ is a square\footnote{We thank Alex Iosevich for pointing out the necessity for $-1$ to not be a square.}, thus $i^2 = -1$ for some $i \in F$, then the set $P:= \{ (x,ix): x \in F \}$ has $\Delta(P) = \{0\}$ even though $|P| = |F|$.  To avoid this degenerate case we assume that $-1$ is not a square, then any two distinct points will have a non-zero distance.  From the fact that any two ``circles'' intersect in at most two points, it is then possible to use extremal graph theory to obtain the
bound
$$ |\Delta(P)| \geq c |P|^{1/2};$$
see also \cite{Erdos-distance}.  This bound is sharp if one takes $P = F^2$, so that $\Delta(P)$ is essentially all of $F$.  Similarly if
one takes $P = G^2$ for any subfield $G$ of $F$.  However, as in the previous section one can hope to improve this bound when
no subfields are available.

From the obvious identity
$$ \Delta(A \times A) = (A - A)^2 + (A - A)^2$$
it is clear that this problem has some connection to the sum-product estimate.  Indeed, any improvement to the trivial bound on $|\Delta(P)|$ can be used (in combination with Lemma \ref{kt} and Lemma \ref{iterated}) to obtain a bound of the form in Lemma \ref{sum-product}.
We now present the converse implication, using the sum-product bounds already obtained to derive a new bound on the distance problem.

\begin{theorem}  Let $F = \Z/p\Z$ for some prime $p$ with $p=3 \mod 4$ (so $-1$ is not a square), and let $P$ be a subset of $F^2$ of cardinality $|P| = N = |F|^\alpha$ for some
$0 < \alpha < 2$.  Then we have
$$ |\Delta(P)| \gtrsim N^{1/2 + \eps}$$
for some $\eps = \eps(\alpha) > 0$.
\end{theorem}

{\bf Remark.}  In the Euclidean analogue to this problem, with $N$ points in $\R^2$, it is conjectured \cite{Erdos-distance} that the above estimate is true for all $\eps < 1/2$.  Currently, this is known for all $\eps < \frac{4e}{5e-1} - \frac{1}{2} \approx 0.364$ \cite{toth}.  However, the Euclidean results depend (among other things) on crossing number technology and thus do not seem to obviously extend to the finite field case.

\begin{proof}
We shall exploit the Szemer\'edi-Trotter-type estimate in Theorem \ref{sz} in much the same way that the actual Szemer\'edi-Trotter theorem \cite{Szemeredi:trotter} was exploited in \cite{cst} for the Euclidean version of Erd\"os's distance problem, or how a Furstenburg set estimate was used in \cite{katz:falc}, \cite{katz:falc-expos} to imply a Falconer
distance set problem result.  The key geometric observation is that the set of points which are equidistant from two fixed points lie on a line (the perpendicular bisector of the two fixed 
points).

We may assume that $|F|$ and $|P|$ are large; in particular, we may assume that $F$ has characteristic greater than 2.  Fix $N$, and suppose for contradiction that
$$ |\Delta(P)| \lesssim N^{1/2+\eps}$$
for some small $0 < \eps \ll 1$ to be chosen later.  For any point $p \in P$, we clearly have the identity
$$\{ (p',r) \in P \times \Delta(P): d(p,p') = r \} | = |P| = N$$
so by Lemma \ref{cz}
$$ \{ (p',p'',r) \in P \times P \times \Delta(P): d(p,p') = d(p,p'') = r; p' \neq p'' \} | \gtrsim N^{3/2 - \eps}.$$
We can of course eliminate the $r$ variable:
$$ \{ (p',p'') \in P \times P: d(p,p') = d(p,p''); p' \neq p'' \} | \gtrsim N^{3/2 - \eps}.$$

Summing this over all $p \in P$ and rearranging, we obtain
$$ \sum_{p',p'' \in P: p' \neq p''} | \{ p \in P: d(p,p') = d(p,p'') \} | \geq c N^{5/2 - \eps}.$$
By the pigeonhole principle, there thus exists $p_0 \in P$ such that
$$ \sum_{p' \in P: p' \neq p_0} | \{ p \in P: d(p,p') = d(p,p_0) \} | \geq c N^{3/2 - \eps}.$$
By translation invariance we may take $p_0 = (0,0)$.  Writing $p' = (a,b)$ and $p = (x,y)$, this becomes
$$ \sum_{(a,b) \in P: (a,b) \neq (0,0)} | \{ (x,y) \in P: (x-a)^2 + (y-b)^2 = x^2 + y^2 \} | \geq c N^{3/2 - \eps}.$$
Thus if we let $l(a,b)$ denote the perpendicular bisector of $(0,0)$ and $(a,b)$:
$$ l(a,b) := \{ (x,y) \in F^2: (x-a)^2 + (y-b)^2 = x^2 + y^2 \} = \{ (x,y) \in F^2: 2ax + 2by = a^2 + b^2 \}$$
and let $L$ be the collection of lines $\{ l(a,b): (a,b) \in P \backslash \{0,0\} \}$, then we have
$$ \{ (p,l) \in P \times L: p \in l \} \geq c N^{3/2 - \eps}.$$
But since all the lines $l(a,b)$ are distinct, we have $|L| = N-1$, while $|P|=N$.  Thus this clearly contradicts Theorem \ref{sz}, and we are done.
\end{proof}

\section{Application to the three-dimensional Kakeya problem}\label{kakeya-sec}

Let $F := \Z/q\Z$ for some prime $q$.  We now use the two-dimensional Szemer\'edi-Trotter theorem to obtain a three-dimensional estimate on Besicovitch sets.

{\bf Definition}  A \emph{Besicovitch set} $P \subseteq F^3$ is a set of points which contains a line in every direction.

The \emph{Kakeya set conjecture for finite fields}\footnote{This is the weakest of the Kakeya conjectures; and corresponds to the ``Minkowski dimension'' form of the Kakeya conjectures in Euclidean space.  There is also a Hausdorff dimension analogue in finite fields, as well as a ``maximal function'' statement; see \cite{gerd:finite} for further discussion.} asserts that for every Besicovitch set $P$, one has the estimate $|P| \geq C_\eps |F|^{3-\eps}$ for every $\eps > 0$; see \cite{wolff:survey}, \cite{gerd:finite} for further discussion on this conjecture.  Previously, the best known lower bound was $|P| \gtrsim |F|^{5/2}$, obtained in \cite{wolff:survey} (see also \cite{wolff:kakeya}, \cite{gerd:finite}).  The purpose of this section is to improve this bound to $|P| \gtrsim |F|^{5/2+\eps}$ for some absolute constant $\eps > 0$.

In fact, we can prove a somewhat stronger statement.  We say that a collection $L$ of lines in $F^3$ obey the \emph{Wolff axiom} if for every 2-plane $\pi$, the number of lines in $l$ which lie in $\pi$ is at most $O(|F|)$.  We will then show

\begin{theorem}\label{kakeya}  Let $L$ be a collection of lines in $F^3$ which obey the Wolff axiom and have cardinality $|L| \sim |F|^2$.  Let $P$ be a collection of points in $F^3$ which contains every line in $L$.  Then $|P| \gtrsim |F|^{5/2+\eps}$ for some absolute constant $\eps > 0$.
\end{theorem}

To see how this theorem implies the claimed bound on Besicovitch sets, observe that a collection of lines consisting of one line in each direction automatically obeys the Wolff axiom and has cardinality $\sim |F|^2$.  

When $\eps = 0$ this bound was obtained in \cite{wolff:survey}, \cite{gerd:finite}.  It was observed in \cite{katzlabatao}, \cite{gerd:finite} that if one replaced the finite field $\Z/q\Z$ by a finite field $F$ which contained a subfield $G$ of index 2, then one could obtain a ``Heisenberg group'' counterexample which showed that Theorem \ref{kakeya} must fail for that field.  Thus, as with the previous results, this theorem must somehow use the non-existence of non-trivial subfields of $F$. It is plausible that one could use the Euclidean analogue \cite{bourgain:ring} of Theorem \ref{sum-product} to prove a similar result in Euclidean space (which would provide a completely different proof of the result in \cite{katzlabatao}), but we do not pursue this question here.

\begin{proof}  Let $0 < \eps \ll 1$ be chosen later.  Fix $L$, $P$, and assume for contradiction that
\begin{equation}\label{p-bound}
 |P| \lesssim |F|^{5/2+\eps}.
\end{equation}

As usual we first begin by running some popularity arguments.  Fix $L$, $P$.  For each point $p \in P$ define the multiplicity $\mu(p)$ by
$$ \mu(p) := | \{ l \in L: p \in l \} |;$$
since every line contains exactly $|F|$ points, we thus have
$$ \sum_{p \in P} \mu(p) = |L| |F| \sim |F|^3.$$
By \eqref{p-bound} the popularity argument, if we thus set $P' \subseteq P$ to be the set of points $p \in P$ where
$$ \mu(p) \gtrsim |F|^{1/2-\eps},$$
then we have
$$ \sum_{p \in P'} \mu(p) \gtrsim |F|^3.$$
For each line $l \in L$ define the multiplicity $\lambda(l)$ by
$$ \lambda(l) := | \{ p \in P': p \in l \}|.$$
Then we can rewrite the previous estimate as
$$ \sum_{l \in L} \lambda(l) \gtrsim |F|^3.$$
Thus if we set $L' \subseteq L$ to be the set of lines $l \in L$ such that
$$ \lambda(l) \gtrsim |F|,$$
then we have
$$ \sum_{l \in L'} \lambda(l) \gtrsim |F|^3.$$
Now define
$$ \mu'(p) := | \{ l \in L': p \in l \} |,$$
so that
$$ \sum_{p \in P'} \mu'(p) \gtrsim |F|^3.$$
Thus by the popularity argument again, if we set $P'' \subseteq P'$ to be the set of points where
$$ \mu'(p) \gtrsim |F|^{1/2-\eps},$$
then
$$ \sum_{p \in P''} \mu'(p) \gtrsim |F|^3.$$
Now define
$$ \lambda'(l) := | \{ p \in P'': p \in l\}|,$$
so that
$$ \sum_{l \in L'} \lambda'(l) \gtrsim |F|^3.$$
If we set $L'' \subseteq L'$ to be set of lines $l$ such that
$$ \lambda'(l) \gtrsim |F|,$$
then we have
$$ \sum_{l \in L''} \lambda'(l) \gtrsim |F|^3,$$
and thus
$$ |\{ (p,l) \in P'' \times L'': p \in l \}| \gtrsim |F|^3.$$
By Lemma \ref{cz} and \eqref{p-bound} we thus have
$$ |\{ (p,l,l') \in P'' \times L'' \times L'': p \in l,l'; l \neq l' \}| \gtrsim |F|^{7/2-\eps}.$$
Define a relation $\sim$ on $L''$ by defining $l \sim l'$ if $l \neq l'$ and the lines $l$ and $l'$ intersect at a point in $P''$.  From the previous we thus have
$$ |\{ (l,l') \in L'' \times L'': l \sim l' \}| \gtrsim |F|^{7/2-\eps},$$
since two lines intersect in at most one point.  Applying Lemma \ref{cz} again with the bound $|L''| = O(|F|^2)$ we obtain
$$ |\{ (l, l_0, l_1) \in L'' \times L'' \times L'': l \sim l_0, l_1; l_0 \neq l_1 \}| \gtrsim |F|^{5-C\eps}.$$
By the pigeonhole principle we can thus find distinct $l_0, l_1 \in L'$ such that
$$ |\{ l \in L'': l \sim l_0, l_1 \}| \gtrsim |F|^{1-C\eps}.$$
Fix $l_0, l_1$, and let $L_* \subseteq L''$ denote the set
$$ L_* := \{ l_* \in L'': l_* \sim l_0, l_1 \},$$
thus we have
\begin{equation}\label{lstar-card} |L_*| \gtrsim |F|^{1-C\eps}.
\end{equation}
Later on we shall complement this lower bound on $L_*$ with an upper bound.

The strategy of the proof will be to pass from $L$ (which is in some sense a two-dimensional subset of a four-dimensional algebraic variety - a Grassmannian, in fact), to $L'$ (which will essentially be a one-dimensional subset of a two-dimensional algebraic variety - namely, the set of lines intersecting both $l_0$ and $l_1$).  The latter situation is much closer to the incidence problem considered in Theorem \ref{sz}, and we will be able to apply that theorem after some algebraic transformations and combinatorial estimates.

For future reference we observe the following non-concentration property of the collection of lines $L'$ (and hence of its subsets $L''$ and $L_*$).

\begin{lemma}\label{regulus}  Let $l^1, l^2, l^3$ be any non-intersecting lines in $F^3$ (not necessarily in $L$).  Then
$$ | \{ l \in L': l^1, l^2, l^3 \hbox{ all intersect } l \} | \lesssim |F|^{1/2+C\eps}.$$
\end{lemma}

As a particular corollary of this lemma, we see that for each 2-plane $\pi$ there are at most $O(|F|^{1/2+C\eps})$ lines in $L'$ which lie in $\pi$.
One further consequence of this is that $l_0$ and $l_1$ are skew (otherwise all the lines in $L_*$ would lie on the plane generated by $l_0$ and $l_1$, and \eqref{lstar-card} would contradict the above corollary).

\begin{proof}  This is a variant of some arguments in \cite{tao:4d} and the
second author.

There are two cases: either some of the lines in $l^1$, $l^2$, $l^3$ are parallel, or they are all mutually skew.  If two of the lines are parallel, ten all the lines $l$ in the above set lie in a plane.  If they are all skew, then it is well known that $l$ lies in a quadratic surface\footnote{A model example is when the lines $l^j$ are of the form $l^j = \{ (x_j, y, x_j y): y \in F \}$ for some distinct $x_1,x_2,x_3 \in F$.  Then all the lines $l$ lie in the quadratic surface $\{ (x,y,xy): x,y \in F \}$.} (i.e. a set of the form $\{ x \in F^3: Q(x) = 0 \}$ for some inhomogeneous quadratic polynomial $Q$) known as the \emph{regulus} generated by $l^1$, $l^2$, $l^3$; see e.g. \cite{schlag:kakeya}, \cite{tao:4d}.  Thus in either case, all the lines $l$ of interest lie inside an algebraic surface $S$ which is either a plane or a quadratic surface.  It will then suffice to show that the set
$$ L_S := \{ l \in L': l \subset S \}$$
has cardinality at most $O(|F|^{1/2+C\eps})$.

We first observe that there are at most $O(|F|)$ lines in $L$ which lie in $S$.  When $S$ is a plane this is just the Wolff axiom.  When $S$ is a quadratic surface this is simply because a quadratic surface contains at most $O(|F|)$ lines.  In particular we have the crude bound $|L_S| = O(|F|)$.

By definition of $L'$ we have
$$ |\{(p,l) \in (P' \cap S) \times L_S: p \in l \}| \gtrsim |F| |L_S|.$$
On the other hand, by Lemma \ref{easy-incidence} we have
$$ |\{ (p,l) \in (P' \cap S) \times L_S: p \in l \}| \lesssim |P' \cap S|^{1/2} |L_S| + |P' \cap S|.$$
Combining the two bounds, we obtain after some algebra
$$ |P' \cap S| \gtrsim \min(|F| |L_S|, |F|^2) \gtrsim |F| |L_S|.$$
By the definition of $P'$, we thus have
$$ |\{(p,l) \in (P' \cap S) \times L: p \in l \}| \gtrsim |F|^{3/2-\eps} |L_S|.$$
The line $l$ certainly intersects $S$, but it need not be contained in $S$.  By the triangle inequality we have either
\begin{equation}\label{in-s} 
|\{(p,l) \in (P' \cap S) \times L: p \in l; l \subset S \}| \gtrsim |F|^{3/2-\eps} |L_S|
\end{equation}
or
\begin{equation}\label{out-s} 
|\{(p,l) \in (P' \cap S) \times L: p \in l; l \not \subset S \}| \gtrsim |F|^{3/2-\eps} |L_S|.
\end{equation}

Suppose first that \eqref{out-s} holds.  Since $S$ is either a plane or a quadratic surface, and $l$ is not contained in $S$, it is clear that $l$ intersects $S$ in at most two places.  Thus the left-hand side of \eqref{out-s} is bounded by at most $2|L| = O(|F|^2)$, and the desired bound $|L_S| = O(|F|^{1/2 + C\eps})$ follows.  Now suppose that \eqref{in-s} holds.
But then there are at most $O(|F|)$ lines $l$ which lie in $S$, and each of those lines $l$ contains at most $|F|$ points $p$.  Thus the left-hand side of \eqref{in-s} is boudned by $O(|F|^2)$, and the desired bound $|L_S| = O(|F|^{1/2 + C\eps})$ again follows.
\end{proof}

Let $l_*$ be any line in $L_*$.  Then $l_0$, $l_1$ each intersect $l_*$ in exactly one point.  Since $l_*$ lies in $L'$, we see from the definition of $L''$ that
$$ | \{ p \in P'': p \in l_*; p \not \in l_0, l_1 \} | \gtrsim |F|.$$
But then by the definition of $P''$, we thus see that
$$ | \{ (p,l) \in P'' \times L': p \in l_*; p \not \in l_0, l_1; p \in l; l \neq l_* \} | \gtrsim |F|^{3/2-C\eps}.$$
Observe that $p$ is uniquely determined by $l$ in the above set.  Define $H(l_*) \subseteq L'$ to be the set of all $l \in L'$ which intersect $l_*$ in a point distinct from where $l_0$ or $l_1$ intersects $l_*$; in the terminology of \cite{wolff:kakeya}, $H(l_*)$ is the \emph{hairbrush} with stem $l_*$.  Then the previous estimate implies that
\begin{equation}\label{hl-lower}
 |H(l_*)| \gtrsim |F|^{3/2-C\eps}
\end{equation}
for all $l_* \in L_*$.

We now complement that lower bound with an upper bound.

\begin{lemma}\label{hb} For each $l_* \in L_*$, we have
$$ |H(l_*)| \lesssim |F|^{3/2+C\eps}.$$
\end{lemma}

\begin{proof}
We use the ``hairbrush'' argument of Wolff \cite{wolff:kakeya}, \cite{wolff:survey}, using the formulation in \cite{gerd:finite}.  All the lines $l$ in $H(l_*)$ intersect $l_*$ but are not coincident to $l_*$.  Thus we have
$$ | \{ (p,l) \in P \times H(l_*): p \in l; p \not \in l_* \} | = |H(l_*)| (|F|-1) \sim |H(l_*)| |F|.$$
By \eqref{p-bound} and Lemma \ref{cz} we thus have
$$ | \{ (p,l,l') \in P \times H(l_*) \times H(l_*): p \in l,l'; p \not \in l_*; l \neq l' \} | \gtrsim |H(l_*)|^2 |F|^{-1/2-\eps}.$$
Now fix $l \in H(l_*)$, and consider how many lines $l'$ could contribute to the above sum.  From the various constraints on $p$, $l$, $l'$, $l_*$ we see that $l$, $l'$, $l_*$ form a triangle, and thus $l'$ lies on the plane generated by $l$ and $l_*$.  By the Wolff axiom there are thus at most $O(|F|)$ choices for $l'$ for any fixed $l$.  Since $l$ and $l'$ clearly determine $p$, we thus have
$$ | \{ (p,l,l') \in P \times H(l_*) \times H(l_*): p \in l,l'; p \not \in l_*; l \neq l' \} | \lesssim |H(l_*)| |F|.$$
Combining this with the previous bound we obtain the Lemma.
\end{proof}

We now refine these hairbrushes slightly.  First we count how many lines $l$ in $H(l_*)$ could intersect $l_0$.  Such lines would lie in the plane generated by $l_*$ and $l_0$ and so there are at most $O(|F|)$ of them by the Wolff axiom (or one could use Lemma \ref{regulus}).  Similarly there are at most $O(|F|)$ lines in $H(l_*)$ which intersect $l_1$.  Thus if we define $\tilde H(l_*)$ to be those lines in $H(l_*)$ which do not intersect either $l_0$ or $l_1$, then we have from \eqref{hl-lower} that
\begin{equation}\label{hl-lower-t}
 |\tilde H(l_*)| \gtrsim |F|^{3/2-C\eps}
\end{equation}

From \eqref{hl-lower-t} we have
$$ |\{ (l_*, l) \in L_* \times L': l \in \tilde H(l_*) \}| \gtrsim |F|^{3/2 - C\eps} |L_*|.$$
By Lemma \ref{cz} (using the bound $|L'| \sim |F|^2$) we thus have
$$ |\{ (l_*, l'_*, l) \in L_* \times L_* \times L': l \in \tilde H(l_*), \tilde H(l'_*); l_* \neq l'_* \}| \gtrsim |F|^{1 - C\eps} |L_*|^2.$$
By the pigeonhole principle we may thus find an $l_* \in L_*$ such that
$$ |\{ (l'_*, l) \in L_* \times L': l \in \tilde H(l_*), \tilde H(l'_*); l_* \neq l'_* \}| \gtrsim |F|^{1 - C\eps} |L_*|.$$
Fix $l_*$.  If we write $L'_* := L_* - \{l_*\}$, we thus have
$$
|\{ (l'_*, l) \in L'_* \times H(l_*): l \in \tilde H(l'_*) \}|\gtrsim |F|^{1 - C\eps} |L_*|.
$$
Let us eliminate some degenerate lines in $L'_*$.  Consider first the contribution of the lines $l'_*$ which contain the point $l_* \cap l_0$.  Such lines also intersect $l_1$ and thus must lie on the plane containing $l_*$ and $l_1$.  The line $l$ must also lie in this plane.  Since this plane contains at most $O(|F|^{1/2+C\eps})$ lines in $L'$ by the corollary to Lemma \ref{regulus}, we thus see that the total contribution of this case is at most $|F|^{1+C\eps}$, which is much smaller than $|F|^{1-C\eps} |L_*|$ by \eqref{lstar-card}.  A similar argument deals with those lines $l'_*$ which contain the point $l_* \cap l_1$.  Thus if we define $L''_*$ to be those lines in $L'_*$ which do not contain either $l_* \cap l_0$ or $l_* \cap l_1$, then we have
$$
|\{ (l'_*, l) \in L''_* \times H(l_*): l \in \tilde H(l'_*) \}|\gtrsim |F|^{1 - C\eps} |L_*|.
$$
If we thus define the quantity
$$ \mu_*(l) := | \{ l'_* \in L''_*: l \in \tilde H(l'_*) \}|$$
then we have
$$ \sum_{l \in \tilde H(l_*)} \mu_*(l) \gtrsim |F|^{1-C\eps} |L_*|.$$
From Lemma \ref{hb} and the popularity argument, if we thus set
$$ H' := \{ l \in \tilde H(l_*): \mu_*(l) \gtrsim |F|^{-1/2-C\eps} |L_*| \}$$
then we have
\begin{equation}\label{h-prime} \sum_{l \in H'} \mu_*(l) \gtrsim |F|^{1-C\eps} |L_*|.
\end{equation}

From Lemma \ref{regulus} we see that 
$$\mu_*(l) \lesssim |F|^{1/2+C\eps}.$$

Indeed, if one unravels all the definitions, we see that the lines $l'_*$ which appear in the definition of $\mu_*(l)$ lie in $L'$ and also intersect the three disjoint lines $l$, $l_0$, $l_1$.

From this bound and \eqref{h-prime} we have a lower bound on $H'$:
\begin{equation}\label{hp-lower}
 |H'| \gtrsim |F|^{1/2-C\eps} |L_*|.
\end{equation}
Comparing this with Lemma \ref{hb} we thus obtain an upper bound on $L_*$:
\begin{equation}\label{lstar-upper}
|L_*| \lesssim |F|^{1 + C\eps}.
\end{equation}

We now perform an algebraic transformation to convert this three-dimensional problem into a two-dimensional problem.  To avoid confusion we shall use the boldface font to denote two-dimensional quantities.

\begin{lemma}\label{2d}  There is a map $\Pi$ from lines $l'_*$ in $L''_*$ to points $\Pi(l'_*)$ in the plane $F^2$, and a map $\Lambda$ from lines $l$ in $H'$ to lines $\Lambda(l)$ in the plane $F^2$, with the following properties:
\begin{itemize}
\item $\Pi$ is injective on $L''_*$.
\item For each line ${\bf l}$ in $F^2$, the fiber $\Lambda^{-1}({\bf l})$ has cardinality at most $|F|^{1/2+C\eps}$.
\item If $(l'_*, l) \in L''_* \times H'$ is such that $l \in \tilde H(l'_*)$, then the point $\Pi(l'_*)$ lies on the line $\Lambda(l)$.
\end{itemize}
\end{lemma}

Let us assume this lemma for the moment, and conclude the proof of the Theorem.
Define the two-dimensional set of points ${\bf P} := \Pi(L''_*)$ and the two-dimensional set of lines ${\bf L} := \Lambda(H')$.  From \eqref{lstar-upper} we have
\begin{equation}\label{p-upper}
 |{\bf P}| \leq |L_*| \lesssim |F|^{1+C\eps}.
\end{equation}
From the multiplicity of $\Lambda$ and \eqref{hp-lower}, \eqref{lstar-card} we have
\begin{equation}\label{l-lower}
|{\bf L}| \gtrsim |F|^{-1/2-C\eps} |H'|
\gtrsim |F|^{-C\eps} |L_*| \gtrsim |F|^{1-C\eps}.
\end{equation}
Now for each $l \in H'$, we have
$$ | \{ l'_* \in L''_*: l \in \tilde H(l'_*) \}| \gtrsim |F|^{-1/2-C\eps} |L_*|
\gtrsim |F|^{1/2-C\eps}$$
by \eqref{lstar-card} and the definitions of $H'$ and $\mu'$.  Applying the first and third parts of Lemma \ref{2d}, we thus conclude that
\begin{equation}\label{heavy}
\{ {\bf p} \in {\bf P}: {\bf p} \in {\bf l} \}| \gtrsim |F|^{1/2-C\eps}
\end{equation}
for all ${\bf l} \in {\bf L}$.

Set $N := |{\bf P}|$.  From \eqref{p-upper}, \eqref{l-lower} we can find a subset ${\bf L'}$ of ${\bf L}$ with 
$$ |F|^{-C\eps} N \lesssim |{\bf L'}| \leq N.$$
Combining this with \eqref{heavy} and then \eqref{p-upper} we obtain
$$ \{ ({\bf p},{\bf l}) \in {\bf P} \times {\bf L'}: {\bf p} \in {\bf l} \}| \gtrsim |F|^{1/2-C\eps} N \gtrsim N^{3/2-C\eps}.$$
But this contradicts the Szemer\'edi-Trotter type estimate in Theorem \ref{sz} if $\eps$ was chosen sufficiently small.  This proves the theorem.

It remains to verify Lemma \ref{2d}.  It is convenient to work in co-ordinates, and for this we shall first normalize the three lines $l_0$, $l_*$, $l_1$.  Recall that $l_0$ and $l_1$ are skew, while $l_*$ intersects both $l_0$ and $l_1$.  After an affine linear transformation, we may set
\begin{align*}
l_0 &:= \{ (x,0,0): x \in F \}\\
l_* &:= \{ (0,0,z): z \in F \}\\
l_1 &:= \{ (0,y,1): y \in F \}
\end{align*}
Consider a line $l'_* \in L''_*$.  This must intersect the line $l_0$ in some point $(x,0,0)$ with $x \neq 0$, and intersect the line $l_1$ in some $(0,y,1)$ with some $y \neq 0$.  Thus $l'_*$ has the form
$$ l'_* = \{ ((1-t)x, ty, t): t \in F \}.$$
We define the map $\Pi: L''_* \to F^2$ by
$$ \Pi(l'_*) := (1/x, 1/y).$$
Clearly $\Pi$ is injective (since two points determine a line).

Now we consider a line $l \in H'$.  This line must intersect $l_*$ at some point $(0,0,z)$ with $z \neq 0,1$.  Thus $l$ takes the form
\begin{equation}\label{l-form}
 l = \{ (a(t-z), b(t-z), t): t \in F \}
\end{equation}
for some $a,b \in F$; note that $a,b \neq 0$ since $l$ is disjoint from $l_0$ and $l_1$.  Suppose this line $l$ intersects the line $l'_*$ mentioned earlier.  Then we must have
$$ a(t-z) = (1-t)x; \quad b(t-z) = ty$$
for some $t \in F$.  Using some algebra to eliminate $t$, we eventually end up with the constraint
$$ xy + (bz-b) x + yaz = 0$$
or (dividing by the non-zero quantity $xy$)
$$ 1 + (bz-b) \frac{1}{y} + az \frac{1}{x} = 0.$$
Thus if we set $\Lambda(l)$ to be the line
$$ \Lambda(l) := \{ (X,Y) \in F^2: 1 + (bz-b) Y + az X = 0 \}$$
(note that this is indeed a line since $a,b \neq 0$ and $z \neq 0,1$), then we see that $\Pi(l'_*) \in \Lambda(l)$ as desired.

It remains to verify the second property of Lemma \ref{2d}.  We consider a line of the form $\{ (X,Y) \in F^2: 1 + \beta Y + \alpha X = 0 \}$ for some $\alpha,\beta \neq 0$ (since these are the only lines in the image of $\Lambda$), and let $L_{\alpha,\beta}$ denote the inverse image of this line under $\Lambda$ in $H'$; thus $L_{\alpha,\beta}$ consists of all the lines in $H'$ of the form
\eqref{l-form}, where $bz-b = \beta$ and $az = \alpha$.  Observe that all such lines must intersect $l_*$, and must also intersect the lines
$\{ (-\alpha, y, 0): y \in F \}$ and $\{ (x, -\beta, 1): x \in F \}$
(indeed, the intersection points are $(-\alpha,-bz,0)$ and $(a(1-z),-\beta,1)$ respectively).  These three lines are disjoint, so by Lemma \ref{regulus} we have $|L_{\alpha,\beta}| \lesssim |F|^{1/2+C\eps}$ as desired.  This concludes the proof of Theorem \ref{kakeya}.
\end{proof}

\end{document}